\documentclass[preprints,article,accept,moreauthors,pdftex]{my-mdpi} 

\firstpage{1} 
\pubvolume{9}
\issuenum{16}
\articlenumber{1883}
\pubyear{2021}
\copyrightyear{2021}
\externaleditor{Academic Editor: Eva Kaslik} 
\datereceived{2 July 2021} 
\daterevised{29 July and 1 August 2021} 
\dateaccepted{6 August 2021} 
\datepublished{8 August 2021} 
\hreflink{https://doi.org/10.3390/math9161883} 
\doinum{10.3390/math9161883}

\pdfoutput=1

\Title{Pontryagin Maximum Principle for Distributed-Order Fractional Systems}

\TitleCitation{Pontryagin Maximum Principle for Distributed-Order Fractional Systems}

\Author{Fa\"{\i}\c{c}al Nda\"{\i}rou $^{\dagger,\ddagger}$\orcidA{} 
and Delfim F. M. Torres *$^{,\ddagger}$\orcidB{}}

\AuthorNames{Fa\"{\i}\c{c}al Nda\"{\i}rou and Delfim F. M. Torres}

\AuthorCitation{Nda\"{\i}rou, F.; Torres, D.F.M.}

\address[1]{Center for Research and Development in Mathematics and Applications (CIDMA),
Department of Mathematics, University of Aveiro, 3810-193 Aveiro, Portugal;
faical@ua.pt}

\corres{\hangafter=1 \hangindent=1.05em \hspace{-0.82em} Correspondence: delfim@ua.pt; Tel.: +351-234-370-668}

\firstnote{\hangafter=1 \hangindent=1.05em \hspace{-0.82em} This 
research is part of first author's Ph.D. project, 
which is carried out at the University of Aveiro under 
the Doctoral Program in Applied Mathematics
of Universities of Minho, Aveiro, and Porto (MAP-PDMA).} 

\secondnote{\hangafter=1 \hangindent=1.05em \hspace{-0.82em} These authors contributed equally to this work.}

\abstract{We consider distributed-order non-local fractional
optimal control problems with controls taking values on a closed set and prove 
a strong necessary optimality condition of Pontryagin type. 
The possibility that admissible controls are
subject to pointwise constraints is new and requires more 
sophisticated techniques to include a maximality condition.
We start by proving results on continuity of solutions due to needle-like control perturbations. 
Then, we derive a differentiability result on the state solutions with respect to the perturbed
trajectories. We end by stating and proving the Pontryagin maximum principle for distributed-order
fractional optimal control problems, illustrating its applicability with an example.}



\allowdisplaybreaks

\begin{document}
	

\section{Introduction}

The idea to consider fractional order systems of distributed order goes back 
to Caputo and the study of anomalous diffusion in viscoelasticity \cite{caputo2}.
The interest on the new operator slowly increased, in particular
with the works of Chechkin et al. \cite{Rev01}, who applied 
distributed order fractional derivatives to study retarding sub-diffusion and
accelerating super-diffusion; Naber \cite{Rev02} studied distributed-order fractional 
subdiffusion processes with different decay rates;
Kochubei \cite{Rev03} applied distributed-order operators 
to the study of ultraslow diffusion; and Mainardi et al. \cite{Rev04} 
applied distributed order fractional derivatives to study Gaussian diffusion. 
The subject is today under strong current research, partially explained
by their relation with physical processes lacking temporal scaling \cite{MR4278073}
and complex non-linear systems \cite{MR4275030}. 
Indeed, the distributed-order definition of the operator 
allows considering superposition of orders and accounting for 
physical phenomena, such as memory effects in composite materials 
and multi-scale effects. A typical example that illustrates 
the capabilities of this class of operators is the mechanical 
behavior of viscoelastic materials
having spatially varying properties.
The literature on experimental applications of
fractional order systems of distributed order is now vast, 
and we refer the interested reader to the review paper of Reference \cite{Rev05}.
For numerical aspects of fractional initial value problems of distributed-order, 
we refer to Reference \cite{MR4276665}.

The calculus of variations is a field of mathematical analysis 
that uses variations, which are small perturbations in functions 
to find maxima and minima of functionals. The Euler--Lagrange equation is 
the main tool for solving such optimization problems, and they have been
developed in the context of fractional calculus to better describe 
non-conservative systems in mechanics \cite{Rev06}.
Necessary optimality condition of Euler--Lagrange type for distributed-order 
problems of the calculus of variations were first introduced and developed in Reference \cite{ricardo}. 
The results were then further generalized by the present authors in Reference \cite{faical},
with the proof of several analytical results and a weak maximum principle of Pontryagin type
for distributed-order fractional optimal control problems. 
Here, we extend and improve 
the theory of optimal control for distributed-order fractional operators initiated in Reference \cite{faical}
by proving a strong version of the Pontryagin maximum principle, which allows the values
of the controls to be constrained to a closed set. The main novelty consists to extend 
the optimality condition proved in Reference \cite{faical} to a maximality condition, 
which yields to the strong version of Pontryagin maximum principle. For this purpose, 
and in contrast with Reference \cite{ricardo,faical}, we use the so-called 
needle-like variations to the control perturbations.

The paper is organized as follows. In Section~\ref{sec:2:P},
we recall some necessary results of the distributed-order fractional calculus.
Our contribution is given in Section~\ref{sectionR}:
we formulate the distributed-order fractional optimal control 
problem under investigation, and we prove the continuity of solutions 
(Lemmas~\ref{conti} and \ref{lem:cor}),
a result on the differentiability of the perturbed trajectories (Lemma~\ref{differenti})
and, finally, the Pontryagin maximum principle (Theorem~\ref{theo}).
We then give an illustrative example of application of the obtained
necessary optimality conditions in Section~\ref{sec:Ex}. 
We end with Section~\ref{sec:conc}, indicating some 
conclusions, the main achievements and novelty of the work,
as well as some future research directions.


\section{Preliminaries}
\label{sec:2:P}

In this section, we recall necessary results and fix notations.
We assume the reader to be familiar with the standard
Riemann--Liouville and Caputo fractional calculi 
\cite{MR3443073,MR1347689}. 

Let $\alpha$ be a real number in $[0, 1]$.
In the sequel, we use the following notation:
\[
L^{\alpha}\left([a, b], \mathbb{R}^n \right)
:=\left\{ x\in L^{1}\left([a, b], \mathbb{R}^n\right)
: I^{\alpha}_{a^{+}} x, I^{\alpha}_{b^{-}} x 
\in AC\left([a, b], \mathbb{R}^n \right) \right\},
\]
where $I^{\alpha}_{a^{+}}$ and $I^{\alpha}_{b^{-}}$ represent, respectively, 
the left and right Riemann--Liouville integral of order $\alpha$.
We also use the notation $AC^{\alpha}\left([a, b], \mathbb{R}^n \right)$ to
represent the set of absolutely continuous functions that can be represented as 
\[
x(t)= x(a) + I^{\alpha}_{a^{+}}f(t) 
\quad \text{ and }\quad x(t)= x(b) + I^{\alpha}_{b^{-}} f(t),
\]
for some functions $f\in L^{\alpha}$.

Let $\psi$ be a non-negative continuous function defined on $[0, 1]$ such that 
$$
\int^1_0 \psi(\alpha)d\alpha >0.
$$ 

This function $\psi$ will act as a distribution of the order 
of differentiation.

\begin{Definition}[See Reference \cite{caputo1}]
The left- and right-sided Riemann--Liouville distributed-order fractional 
derivatives of a function $x\in L^{\alpha}$ are defined, respectively, by 
\[
\mathbb{D}^{\psi(\cdot)}_{a^{+}}x(t)= \int^1_0 \psi(\alpha)
\cdot D^{\alpha}_{a^{+}}x(t)d\alpha 
\quad \text{ and } \quad 
\mathbb{D}^{\psi(\cdot)}_{b^{-}}x(t)
= \int^1_0 \psi(\alpha)\cdot D^{\alpha}_{b^{-}}x(t)d\alpha,
\]
where $D^{\alpha}_{a^{+}}$ and $D^{\alpha}_{b^{-}}$ are, 
respectively, the left- and right-sided Riemann--Liouville 
fractional derivatives of order $\alpha$.
\end{Definition}

\begin{Definition}[See Reference \cite{caputo1}]
The left- and right-sided Caputo distributed-order fractional 
derivatives of a function $x\in AC^{\alpha}$ 
are defined, respectively, by
\[
^{C}\mathbb{D}^{\psi(\cdot)}_{a^{+}}x(t)= \int^1_0 \psi(\alpha)
\cdot ^{C}D^{\alpha}_{a^{+}}x(t)d\alpha 
\quad \text{ and } \quad 
^{C}\mathbb{D}^{\psi(\cdot)}_{b^{-}}x(t)
= \int^1_0 \psi(\alpha)\cdot ^{C}D^{\alpha}_{b^{-}}x(t)d\alpha,
\]
where $^{C}D^{\alpha}_{a^{+}}$ and $^{C}D^{\alpha}_{b^{-}}$ are,
respectively, the left- and right-sided Caputo fractional 
derivatives of order $\alpha$.
\end{Definition}

As noted in Reference \cite{ricardo}, there is a relation between the 
Riemann--Liouville and the Caputo distributed-order fractional derivatives:
\[
^{C}\mathbb{D}^{\psi(\cdot)}_{a^{+}}x(t)
= \mathbb{D}^{\psi(\cdot)}_{a^{+}}x(t)
- x(a)\int^1_0 \frac{\psi(\alpha)}{\Gamma(1-\alpha)}(t-a)^{-\alpha}d\alpha
\]
and 
\[
^{C}\mathbb{D}^{\psi(\cdot)}_{b^{-}}x(t)
= \mathbb{D}^{\psi(\cdot)}_{b^{-}}x(t)- x(b)
\int^1_0 \frac{\psi(\alpha)}{\Gamma(1-\alpha)}(b-t)^{-\alpha}d\alpha.
\]

Along the text, we use the notation
\[
\mathbb{I}^{1-\psi(\cdot)}_{b^{-}}x(t)
=\int^1_0\psi(\alpha)\cdot I^{1-\alpha}_{b^{-}}x(t)d\alpha,
\]
where $I^{1-\alpha}_{b^{-}}$ represents the right 
Riemann--Liouville fractional integral of order $1-\alpha$.

The following results will be useful for our purposes.
In concrete, integration by parts will be used in the
proof of the Pontryagin maximum principle (Theorem~\ref{theo}).

\begin{Lemma}[Integration by parts formula \cite{ricardo}]
\label{lemma1}
Let $x\in L^{\alpha}$ and $y\in AC^{\alpha}$. Then, 
\[
\int^b_a x(t)\cdot ^{C}\mathbb{D}^{\psi(\cdot)}_{a^{+}}y(t)dt 
= \left[ y(t)\cdot \mathbb{I}^{1-\psi (\cdot)}_{b^{-}}x(t) \right]^b_a 
+ \int^b_a y(t)\cdot \mathbb{D}^{\psi(\cdot)}_{b^{-}}x (t)dt.
\]
\end{Lemma}

It follows a generalized Gr\"{o}nwall inequality
that will be used in Section~\ref{sub:sec:CS}.

\begin{Lemma}[Gr\"{o}nwall inequality \cite{gronwall}]
\label{thm:gronwall}
Let $\alpha$ be a positive real number and let $a(\cdot)$, 
$b(\cdot)$, and $u(\cdot)$ be non-negative continuous functions 
on $[0, T]$ with $b(\cdot)$ monotonic increasing on $[0, T)$.
If
\[
u(t)\leq a(t) + b(t)\int^t_0(t-s)^{\alpha-1}a(s)ds,
\]
then 
\[
u(t)\leq a(t) + \int^t_0 \left[ \sum^{\infty}_{n=0}
\frac{\left(b(t)\Gamma(\alpha ) \right)^n}{\Gamma(n\alpha)} 
(t-s)^{n\alpha-1}\right]ds
\]
for all $t\in [0,T)$.
\end{Lemma}


\section{Main Results}
\label{sectionR}

In this work, we look for an essentially bounded control 
$u \in L^{\infty}\left([a, b], \mathbb{R}^m \right)$ 
and the corresponding state trajectory 
$x\in AC^{\alpha}\left([a, b], \mathbb{R}^n \right)$, 
solution to the following distributed-order non-local 
fractional-order optimal~control~problem:
\begin{equation}
\label{pmp}
\begin{gathered}
J[x(\cdot), u(\cdot)]= \int^{b}_{a} L\left(t, x(t), u(t)\right)dt 
\longrightarrow \max,\\
^{C}\mathbb{D}^{\psi(\cdot)}_{{a}^{+}}x(t)
= f\left(t, x(t), u(t)\right), \quad t\in [a, b] \ a.e.,\\
x(\cdot) \in AC^{\alpha}, \quad u(\cdot) \in L^{\infty},\\
x(a)= x_a\in \mathbb{R}^n, \quad u(t)\in \Omega,
\end{gathered}
\end{equation}
where $\Omega$ is a closed subset of $\mathbb{R}^m$. 
The data functions $L: [a, b]\times \mathbb{R}^n \times \mathbb{R}^m \rightarrow \mathbb{R}$ 
and $f: [a, b]\times \mathbb{R}^n \times \mathbb{R}^m \rightarrow \mathbb{R}^n$ 
are subject to the following assumptions:
\begin{itemize}
\item The function $f$ is continuous in all its three arguments.
\item The function $f$ is continuously differentiable with respect to state variable $x$ and, 
in particular, locally Lipschitz-continuous, that is, for every compact 
$B\subset \mathbb{R}^n$ and for all $x, y\in B$ there is $K> 0$ such that
\[
\parallel f(t, x, u)-f(t, y, u) \parallel \leq K \parallel x-y \parallel.
\]
\item With respect to the control $u$, there exists $M>0$ such that 
\[
\parallel f(t, x, u)\parallel \leq M
\quad \forall (t, x)\in [a, b]\times \mathbb{R}^n.
\]
\item The cost integrand $L$ satisfies the same assumptions as $f$.
\end{itemize}


\subsection{Sensitivity Analysis} 
\label{sub:sec:CS}

Now, our concern is to establish continuity and differentiability 
results on the state solutions for any control perturbation 
(Lemmas~\ref{conti}--\ref{differenti}), which 
are then used in Section~\ref{subsec:PMP} to prove a
necessary optimality condition for the optimal control
problem \eqref{pmp}. With this purpose, let us denote by 
$\mathcal{L}\left[ F(\cdot)\right]$ the set of all Lebesgue 
points in $[a, b)$ of the essentially bounded functions 
$t\mapsto f(t, x(t), u(t))$ and $t\mapsto L(t, x(t), u(t))$. Thus, 
let $(\tau, v) \in \mathcal{L}\left[ F(\cdot)\right] \times \Omega$, 
and, for every $\theta \in [0, b-\tau)$, let us consider the 
needle-like variation $u^{\theta} \in L^{\infty}\left([a, b], \mathbb{R}^n \right)$ 
associated to the optimal control $u^{*}$, which is given by
\begin{equation}
u^{\theta}(t)= 
\begin{cases}
u^{*}(t) \quad \text{ if } \quad t\not\in [\tau-\theta, \tau),\\
v \qquad \quad \text{ if } \quad t\in [\tau-\theta, \tau),
\end{cases}
\end{equation}
for almost every $t\in [a, b]$.

\begin{Lemma}[Continuity of solutions]
\label{conti} 
For any $(\tau, v) \in \mathcal{L}\left[ F(\cdot)\right] \times \Omega$, 
denote by $x^{\theta}$ the corresponding state trajectory to the needle-like 
variation $u^{\theta}$, that is, the state solution of 
\[
^{C}\mathbb{D}^{\psi(\cdot)}_{{a}^{+}} x^{\theta}(t)
= f\left(t, x^{\epsilon}(t), u^{\theta}(t)\right), \quad x^{\theta}(a)= x_a.
\]

Then, we have that $x^{\theta}$ converges uniformly to the optimal state 
trajectory $x^{*}$ whenever $\theta$ tends to zero.
\end{Lemma}

\begin{proof}
We have that
\[
^{C}\mathbb{D}^{\psi(\cdot)}_{{a}^{+}}\left( x^{\theta}(t) 
- x^{*}(t)\right) 
=  f(t, x^{\theta}(t), u^{\theta}(t)) 
- f\left(t, x^{*}(t), u^{*}(t)\right).
\]

Then, by definition of the distributed-order operator,
\[
\int^1_0 \psi(\alpha)^CD^{\alpha}_{{a}^{+}}\left( x^{\theta}(t)
- x^{*}(t)\right) d\alpha 
=  f(t, x^{\theta}(t), u^{\theta}(t)) 
- f\left(t, x^{*}(t), u^{*}(t)\right).
\]

Now, using the mean value theorem for integrals, 
there exists an $\bar{\alpha}$ such that 
\[
^CD^{\bar{\alpha}}_{{a}^{+}}\left( x^{\theta}(t)
- x^{*}(t)\right)
= \frac{1}{m}\left[ f(t, x^{\theta}(t), u^{\theta}(t)) 
- f\left(t, x^{*}(t), u^{*}(t)\right)\right]
\]
with 
$$
m= \int^1_0 \psi(\alpha)d\alpha.
$$

Therefore, by the left inverse property, 
we obtain the following integral representation: 
\[
x^{\theta}(t)-x^{*}(t)= \frac{1}{m} I^{\bar{\alpha}}_{a^{+}}\left( 
f(t, x^{\theta}(t), u^{\theta}(t)) 
- f\left(t, x^{*}(t), u^{*}(t)\right)\right).
\]

Moreover, note that 
\begin{multline*}
f(t, x^{\theta}(t), u^{\theta}(t)) 
- f\left(t, x^{*}(t), u^{*}(t)\right)
= \{ f(t, x^{\theta}(t), u^{\theta}(t)) - f(t, x^{*}(t), u^{\theta}(t))\} \\
+ \{f(t, x^{*}(t), u^{\theta}(t)) - f\left(t, x^{*}(t), u^{*}(t)\right)\}.
\end{multline*}

With the help of the triangular inequality, we can write that
\begin{multline*}
\parallel x^{\theta}(t)-x^{*}(t) \parallel \leq \frac{1}{m} 
I^{\bar{\alpha}}_{a^{+}}\left( \parallel f(t, x^{\theta}(t), u^{\theta}(t)) 
-f(t, x^{*}(t), u^{\theta}(t)) \parallel \right)\\
+ \frac{1}{m} I^{\bar{\alpha}}_{(\tau-\theta)^{+}}\left( 
\parallel f(t, x^{*}(t), u^{\theta}(t))
- f\left(t, x^{*}(t), u^{*}(t)\right) \parallel \right),
\end{multline*}
since $u^{\theta}$ and $u^{*}$ are different only 
on $[\tau-\theta, \tau]$. From the Lipschitz property of $f$ 
and the boundedness with respect to the control, it follows that
\[
\parallel x^{\theta}(t)-x^{*}(t) \parallel \leq \frac{K}{m} 
I^{\bar{\alpha}}_{a^{+}}\left( \parallel x^{\theta}(t)
-x^{*}(t) \parallel \right) + \frac{1}{m}
\cdot 2M\frac{\theta^{\bar{\alpha}}}{\Gamma(\alpha + 1)}.
\]

Now, by applying the fractional Gr\"{o}nwall inequality
(Lemma~\ref{thm:gronwall}), it follows that
\begin{align*}
\parallel x^{\theta}(t)-x^{*}(t) \parallel 
\leq \frac{2M\theta^{\bar{\alpha}}}{m\Gamma(\alpha + 1)}\left[   
1+ \int_a^t \sum^{\infty}_{n=1}\frac{K^n}{\Gamma(n\bar{\alpha})} 
(t-s)^{n\bar{\alpha}-1}ds \right] \leq \varpi_1 \theta^{\bar{\alpha}},
\end{align*}
where $\displaystyle{ \varpi_1 = \frac{2M}{m\Gamma(\alpha + 1)}
E_{\alpha, 1}\left( K(b-a)^{\alpha}\right) }$, and $E_{\alpha, 1}$ 
is the Mittag--Leffler function of parameter $\bar{\alpha}$. Hence, 
by taking the limit when $\theta$ tends to zero, we obtain the desired result: 
$x^{\theta} \rightarrow x^{*}$ for all $t\in [a, b]$.
\end{proof}

The next result is a corollary of Lemma~\ref{conti}.

\begin{Lemma}
\label{lem:cor}
There exists $\varpi_2 \geq 0$ such that 
\[
\parallel x^{\theta}(t)-x^{*}(t) \parallel \leq \varpi_2 \theta
\left(t-(\tau-\theta)\right)^{\bar{\alpha}-1} 
\quad \forall t\in \, ]\tau-\theta, b].
\]
\end{Lemma}

\begin{proof}
Using similar arguments of Lipschitz-continuity of $f$ 
and its boundedness with respect to the control $u$, we get
\begin{multline*}
\parallel x^{\theta}(t)-x^{*}(t) \parallel 
\leq \frac{M}{m\Gamma(\bar{\alpha})}\int^{\tau}_{\tau-\theta}(t-s)^{\bar{\alpha}-1}ds\\ 
+ \frac{K}{m\Gamma(\bar{\alpha})}\int^{\tau}_{\tau-\theta}(t-s)^{\bar{\alpha}-1}
\parallel x^{\theta}(s)-x^{*}(s) \parallel ds
+ \frac{K}{m\Gamma(\bar{\alpha})}\int^{t}_{\tau}(t-s)^{\bar{\alpha}-1}
\parallel x^{\theta}(s)-x^{*}(s) \parallel ds.
\end{multline*}

Note that $\displaystyle{\int^{\tau}_{\tau-\theta}(t-s)^{\bar{\alpha}-1}ds 
\leq \theta \left(t-(\tau-\theta)\right)^{\bar{\alpha}-1} }$, and, 
as a consequence of Lemma~\ref{conti}, we obtain that
\begin{multline*}
\parallel x^{\theta}(t)-x^{*}(t) \parallel \leq \frac{M}{m\Gamma(\bar{\alpha})}(M
+ \varpi_1 K\theta^{\bar{\alpha}})\theta\left(t-(\tau-\theta)\right)^{\bar{\alpha}-1}\\ 
+  \frac{K}{m\Gamma(\bar{\alpha})}\int^{t}_{\tau}(t-s)^{\bar{\alpha}-1}
\parallel x^{\theta}(s)-x^{*}(s) \parallel ds.
\end{multline*}

We conclude the proof by applying again the fractional Gr\"{o}nwall inequality 
(\mbox{Lemma~\ref{thm:gronwall}}), in which we set 
$\varpi_2= \frac{1}{m}M+\varpi_1 K \theta^{\bar{\alpha}}
E_{\alpha, 1}\left( K(b-a)^{\alpha}\right)$.
\end{proof}

\begin{Lemma}[Differentiability of the perturbed trajectory]
\label{differenti}
For all $(\tau, v) \in \mathcal{L}\left[ F(\cdot)\right] \times \Omega$, 
we have that the variational trajectory 
$\displaystyle{\frac{x^{\theta}(\cdot)-x^{*}(\cdot)}{\theta}}$
is uniformly convergent to $\eta(\cdot)$ when $\theta$ tends to zero, 
where $\eta(\cdot)$ is the unique solution to the distributed-order 
left Caputo fractional Cauchy problem
\begin{equation}
\label{eqLinear}
\begin{cases}
\displaystyle{ ^{C}\mathbb{D}^{\psi(\cdot)}_{{\tau}^{+}} \eta(t)
=  \frac{\partial f(t, x^{*}(t), u^{*}(t))}{\partial x} \cdot \eta(t)}, 
\quad t\in ]\tau, b],\\[3mm]
\displaystyle{I^{1-\bar{\alpha}}_{{\tau}^{+}} \eta(\tau)
= \frac{1}{m}\left[f(\tau, x^{*}(\tau), v)
-f(\tau, x^{*}(\tau), u^{*}(\tau))\right]}.
\end{cases}
\end{equation}
\end{Lemma} 

\begin{proof}
Set $\displaystyle{z^{\theta}(t)= \frac{x^{\theta}(t)-x^{*}(t)}{\theta}-\eta(t)}$ 
for all $t\in [\tau, b]$. Our aim is to prove that $z^{\theta}$ 
converges uniformly to zero on $[\tau, b]$ whenever $\theta \rightarrow 0$. 
The integral representation of $z^{\theta}$ is given as follows:
\begin{multline}
\label{eq:Intrepres}
z^{\theta}(t)= -\frac{1}{m\Gamma(\bar{\alpha})}(t-\tau)^{\bar{\alpha}-1}\left( 
f(\tau, x^{*}(\tau), v)-f(\tau, x^{*}(\tau), u^{*}(\tau))\right)\\
+\frac{1}{m\Gamma(\bar{\alpha})}\int^t_{\tau^{+}}(t-s)^{\bar{\alpha}-1}\left[
\frac{f(s, x^{\theta}(s), u^{*}(s))-f(s, x^{*}(s), u^{*}(s))}{\theta} \right.\\
\left. -\frac{\partial f(s, x^{*}(s), u^{*}(s))}{\partial x}
\times \frac{x^{\theta}(s)-x^{*}(s)}{\theta}\right] ds\\
+\frac{1}{m\Gamma(\bar{\alpha})}\int^t_{\tau^{+}}(t-s)^{\bar{\alpha}-1}
\frac{\partial f(s, x^{*}(s), u^{*}(s))}{\partial x}\times z^{\theta}(s) ds
\end{multline}
for every $t\in [\tau, b]$. Let us investigate the two first terms 
of the right-hand side of \eqref{eq:Intrepres}. 
By boundedness of $f$ with respect to $u$, we have that 
\[
\left\Vert-\frac{1}{m\Gamma(\bar{\alpha})}(t
-\tau)^{\bar{\alpha}-1}\left( f(\tau, x^{*}(\tau), v)-f(\tau, x^{*}(\tau),
u^{*}(\tau))\right)\right\Vert 
\leq \frac{2M}{\Gamma(\bar{\alpha})}(b-\tau)^{\bar{\alpha}}.
\]

Further, using the classical Taylor formula with integral rest, we have 
\begin{multline*}
\frac{f(s, x^{\theta}(s), u^{*}(s))-f(s, x^{*}(s), u^{*}(s))}{\theta} 
-\frac{\partial f(s, x^{*}(s), u^{*}(s))}{\partial x}\times 
\frac{x^{\theta}(s)-x^{*}(s)}{\theta}\\
=\int^1_0 \left( \frac{\partial f(s, x^{*}(s)
+w(x^{\theta}(s)-x^{*}(s)), u^{*}(s))}{\partial x}
-\frac{\partial f(s, x^{*}(s), u^{*}(s))}{\partial x}\right)\\
\times \left( \frac{x^{\theta}(s)-x^{*}(s)}{\theta} \right) dw.
\end{multline*}

Hence, from Lemma~\ref{lem:cor}, we deduce that 
$\displaystyle{\left\Vert \frac{x^{\theta}(s)-x^{*}(s)}{\theta} \right\Vert
\leq \varpi_2  \left(t-(\tau-\theta)\right)^{\bar{\alpha}-1}}$. 
Next, we set
\[
\varsigma_{\theta}(s)= \int^1_0 \left\Vert \frac{\partial 
f(s, x^{*}(s)+w(x^{\theta}(s)-x^{*}(s)), u^{*}(s))}{\partial x}
-\frac{\partial f(s, x^{*}(s), u^{*}(s))}{\partial x}\right\Vert ds,
\]
and, referring to Lemma A.3 in Reference \cite{bourdin}, we get an estimate 
for the second term of \eqref{eq:Integral}, and we end the proof 
by application of the fractional Gr\"{o}nwall inequality 
of Lemma~\ref{thm:gronwall}.
\end{proof}


\subsection{Pontryagin's Maximum Principle of Distributed-Order}
\label{subsec:PMP}

It follows the main result of our work:
a distributed-order Pontryagin maximum principle
for the fractional-order optimal control problem \eqref{pmp}.

\begin{Theorem}[Pontryagin Maximum Principle for \eqref{pmp}]
\label{theo}
If $(x^{*}(\cdot), u^{*}(\cdot))$ is an optimal pair for \eqref{pmp}, 
then there exists $\lambda \in L^{\alpha}$, called the adjoint function variable, 
such that the following conditions hold for all $t$ in the interval $[a, b]$:
\begin{itemize}
\item the maximality condition
\begin{equation}
\label{maxcondition}
H(t, x^{*}(t), u^{*}(t), \lambda(t))
= \underset{\omega \in \Omega} \max \, H\left(t, x^{*}(t), \omega, \lambda(t)\right);
\end{equation}
\item the adjoint system
\begin{equation}
\label{adj}
\mathbb{D}^{\psi(\cdot)}_{b^{-}} \lambda(t)
= \frac{\partial H}{\partial x}(t,x^{*}(t), u^{*}(t), \lambda(t));
\end{equation}
\item the transversality condition
\begin{equation}
\label{trans}
\mathbb{I}^{1-\psi(\cdot)}_{b^{-}}\lambda(b)=0,
\end{equation}
\end{itemize}
where the Hamiltonian $H$ is defined by
\[
H(t, x, u, \lambda)= L(t, x, u) + \lambda \cdot f(t, x, u).
\]
\end{Theorem}

\begin{proof}
First of all, note that the regularity of function $f$ with respect to the state variable 
(recall that $f$ is continuously differentiable with respect to $x$) is exactly 
as in our previous paper \cite{faical}. For this reason, the adjoint system \eqref{adj}
and its transversality condition \eqref{trans} remain exactly the same
as the ones proved in Reference \cite{faical}. Therefore, we only need to prove 
the maximality condition \eqref{maxcondition}, which is new due to less 
regularity of $f$ with respect to control functions and the fact that now
the controls take values on the closed $\Omega$ set.
We start by using integration by parts 
(Lemma~\ref{lemma1}) for functions $\lambda \in L^{\alpha}$ 
and $\eta \in AC^{\alpha}$ on $[\tau, b]$:
\begin{equation}
\label{integradjoint}
\int^b_{\tau} \lambda(s)\cdot ^{C}\mathbb{D}^{\psi(\cdot)}_{\tau^{+}}\eta(s)ds 
= \left[ \eta(s)\cdot \mathbb{I}^{1-\psi (\cdot)}_{b^{-}}\lambda(s) \right]^b_{\tau} 
+ \int^b_{\tau} \eta(s)\cdot \mathbb{D}^{\psi(\cdot)}_{b^{-}} \lambda(s)ds,
\end{equation}
where $\lambda$ is the adjoint variable given in Reference \cite{faical}:
\begin{equation}
\label{Adjsystem}
\begin{cases}
\mathbb{D}^{\psi(\cdot)}_{b^{-}} \lambda(t)
= \frac{\partial L}{\partial x}(t,x^{*}(t), u^{*}(t)) 
+ \lambda(t)\cdot \frac{\partial f}{\partial x}(t,x^{*}(t), u^{*}(t)),\\[3mm]
\mathbb{I}^{1-\psi(\cdot)}_{b^{-}}\lambda(b)=0.
\end{cases}
\end{equation}

Substituting \eqref{Adjsystem} and the variational differential system 
given in \eqref{eqLinear} into \eqref{integradjoint}, we obtain~that
\begin{multline*}
\int^b_{\tau} \lambda(s)\cdot\left(\frac{\partial f(s, x^{*}(s), u^{*}(s))}{\partial x}
\cdot \eta(s)\right)ds = -\eta(\tau)\mathbb{I}^{1-\psi(\cdot)}_{b^{-}}\lambda(\tau) \\
+ \int^b_{\tau}\eta(s)\left(\frac{\partial L}{\partial x}(s,x^{*}(s), u^{*}(s)) 
+  \lambda(s)\cdot \frac{\partial f}{\partial x}(s,x^{*}(s), u^{*}(s)) \right)ds,
\end{multline*}
which leads to
\begin{equation}
\label{Lamrelatio}
\eta(\tau)\mathbb{I}^{1-\psi(\cdot)}_{b^{-}}\lambda(\tau) 
= \int^b_{\tau} \eta(s)\cdot \frac{\partial L(s, x^{*}(s), u^{*}(s))}{\partial x}ds.
\end{equation}

Next, recall that, from the definition of distributed-order fractional integral 
and the mean value theorem, we have the existence of an $\bar{\alpha}$ such that\vspace{-6pt}

\begin{equation}
\label{eq:Integral}
\mathbb{I}^{1-\psi(\cdot)}_{b^{-}}\lambda(\tau) 
= \int^1_0\psi(\alpha)I^{1-\bar{\alpha}}_{b^{-}}\lambda(\tau)d\alpha
= mI^{1-\bar{\alpha}}_{b^{-}}\lambda(\tau),
\end{equation}
where $m=\int^1_0 \psi(\alpha)d\alpha$. Moreover, by the fundamental law of calculus 
and the duality of the Riemann-Liouville integral operator, we have also that
\begin{equation*}
\begin{split}
\eta(\tau)I^{1-\bar{\alpha}}_{b^{-}}\lambda(\tau)
&= \frac{d}{d\tau}\left(-\int^b_{\tau}\eta(s)I^{1-\bar{\alpha}}_{b^{-}}\lambda(s)ds\right)\\
&= \frac{d}{d\tau}\left(-\int^b_{\tau}\lambda(s)I^{1-\bar{\alpha}}_{\tau^{+}}\eta(s)ds\right)\\
&=\lambda(\tau)I^{1-\bar{\alpha}}_{\tau^{+}}\eta(\tau).
\end{split}
\end{equation*}

Next, using the boundary condition from system \eqref{eqLinear}, it yields
\begin{equation*}
\begin{split}
\eta(\tau)\mathbb{I}^{1-\psi(\cdot)}_{b^{-}}\lambda(\tau)
&= m\eta(\tau)I^{1-\bar{\alpha}}_{b^{-}}\lambda(\tau)
=m\lambda(\tau)I^{1-\bar{\alpha}}_{\tau^{+}}\eta(\tau)\\
&=m\lambda(\tau)\left(\frac{1}{m}[f(\tau, x^{*}(\tau), v)
-f(\tau, x^{*}(\tau), u^{*}(\tau))] \right),
\end{split}
\end{equation*}
that is, $\eta(\tau)\mathbb{I}^{1-\psi(\cdot)}_{b^{-}}\lambda(\tau)
= \lambda(\tau)\cdot \left(f(\tau, x^{*}(\tau), v)-f(\tau, x^{*}(\tau), u^{*}(\tau)) \right)$.
Finally, substituting this expression into \eqref{Lamrelatio}, we get
\begin{equation}
\label{eq:333}
\lambda(\tau)\cdot \left(f(\tau, x^{*}(\tau), v)-f(\tau, x^{*}(\tau), u^{*}(\tau)) \right)
= \int^b_{\tau} \eta(s)\cdot \frac{\partial L(s, x^{*}(s), u^{*}(s))}{\partial x}ds.
\end{equation}

However, with respect to the cost functional $J$, the limit
\begin{equation}
\label{Jnegative}
\lim_{\theta \rightarrow 0^{+}} \frac{J\left[x^{\theta}(\cdot), u^{\theta}(\cdot) \right]
-J\left[x^{*}(\cdot), u^{*}(\cdot) \right]}{\theta}\leq 0
\end{equation}
because, by assumption, $(x^{*}, u^{*})$ is an optimal pair. This limit can be written as 
\begin{multline}
\label{eq:limit}
\lim_{\theta \rightarrow 0^{+}} \frac{J\left[x^{\theta}(\cdot), u^{\theta}(\cdot) \right]
-J\left[x^{*}(\cdot), u^{*}(\cdot) \right]}{\theta}\\
= \lim_{\theta \rightarrow 0^{+}}\frac{1}{\theta}\int^{\tau}_{\tau -\theta} \left[ 
L(s, x^{*}(s), v)-L(s, x^{*}(s), u ^{*}(s)) \right]ds\\ 
+ \lim_{\theta \rightarrow 0^{+}}\int^b_{\tau} \frac{L(s, x^{\theta}(s), u ^{\theta}(s))
-L(s, x^{*}(s), u ^{\theta}(s))}{\theta}ds.
\end{multline}

Considering the fact that $\tau$ is a Lebesgue point of 
$$
L(s, x^{*}(s), v)-L(s, x^{*}(s), u ^{*}(s)) := \psi(s), 
$$
it follows from the Lebesgue differentiation property
\[
\left| \frac{1}{\theta}\int^{\tau}_{\tau -\theta}\psi(s)ds -\psi(\tau)\right| 
= \left| \frac{1}{\theta}\int^{\tau}_{\tau -\theta}\left(\psi(s)-\psi(\tau)\right)ds \right| 
\leq \frac{1}{\theta}\int^{\tau}_{\tau -\theta}\left|\psi(s)-\psi(\tau)\right|ds
\]
that
\begin{equation}
\begin{split}
\lim_{\theta \rightarrow 0^{+}}\frac{1}{\theta}\int^{\tau}_{\tau -\theta} \left[ 
L(s, x^{*}(s), v)-L(s, x^{*}(s), u ^{*}(s)) \right]ds\\
=L(\tau, x^{*}(\tau), v)-L(\tau, x^{*}(\tau), u ^{*}(\tau)).
\end{split}
\end{equation}
Moreover, with respect to the third limit in \eqref{eq:limit}, 
we can apply the Lipschitz property of $L$ to obtain 
\[
\left|\frac{L(s, x^{\theta}(s), u ^{\theta}(s))-L(s, x^{*}(s), u ^{\theta}(s))}{\theta}\right| 
\leq K \left\Vert \frac{x^{\theta}-x^{*}}{\theta}\right\Vert.
\]

Therefore, because $\displaystyle{ \frac{x^{\theta}-x^{*}}{\theta}}$ is a uniformly convergent 
series of functions, we conclude that the integrand 
$$
\displaystyle{\frac{L(s, x^{\theta}(s), u ^{\theta}(s))-L(s, x^{*}(s), u ^{\theta}(s))}{\theta}}
$$ 
is uniformly bounded. Furthermore, we have 
\begin{multline*}
L(s, x^{\theta}(s), u^{\theta}(t))
= L\left(s, x^{*}(s), u^{\theta}(s)\right) \\
+ (x^{\theta}(s)-x^{*}(s))
\cdot \frac{\partial L(s, x^{*}(s), u^{\theta}(s))}{\partial x} 
+ o\left( \|x^{\theta}-x^{*}\| \right).
\end{multline*}

Next, by the continuity Lemma~\ref{conti}, we have $\|x^{\theta}-x^{*}\| \rightarrow 0$ 
whenever $\theta \rightarrow 0$. Thus, we can express the residue term 
only as function of $\theta$, that is,
\[
L(s, x^{\theta}(s), u^{\theta}(s))
= L\left(s, x^{*}(s), u^{\theta}(s)\right) + (x^{\theta}(s)-x^{*}(s))
\cdot \frac{\partial L(s, x^{*}(s), u^{\theta}(s))}{\partial x} + o(\theta),
\]
and the following expression yields for the second limit:
\begin{equation*}
\begin{split}
\lim_{\theta \rightarrow 0}\frac{L(s, x^{\theta}(s), u^{*}(s))
- L\left(t, x^{*}(t), u^{*}(t)\right)}{\theta} 
&=\frac{\partial L(s, x^{*}(s), u^{*}(s))}{\partial x} 
\cdot \lim_{\theta \rightarrow 0}\frac{(x^{\theta}(s)-x^{*}(s))}{\theta}\\
&=\frac{\partial L(s, x^{*}(s), u^{*}(s))}{\partial x} \cdot \eta(s).
\end{split}
\end{equation*}

Hence, thanks to the Lebesgue bounded convergence theorem, 
\[
\lim_{\theta \rightarrow 0^{+}}\int^b_{\tau} \frac{L(s, x^{\theta}(s), u ^{\theta}(s))
-L(s, x^{*}(s), u ^{\theta}(s))}{\theta}ds
=\frac{\partial L(s, x^{*}(s), u^{*}(s))}{\partial x} \cdot \eta(s)
\]
and, altogether, we get
\begin{multline*}
\lim_{\theta \rightarrow 0^{+}} 
\frac{J\left[x^{\theta}(\cdot), u^{\theta}(\cdot) \right]
-J\left[x^{*}(\cdot), u^{*}(\cdot) \right]}{\theta}\\
= L(\tau, x^{*}(\tau), v)
-L(\tau, x^{*}(\tau), u ^{*}(\tau))
+ \frac{\partial L(s, x^{*}(s), u^{*}(s))}{\partial x} \cdot \eta(s).
\end{multline*}

Hence, using inequality \eqref{Jnegative} and \eqref{eq:333}, we obtain that
\[
L(\tau, x^{*}(\tau), v)-L(\tau, x^{*}(\tau), u ^{*}(\tau)) 
+\lambda(\tau)\cdot \left[f(\tau, x^{*}(\tau), v)
-f(\tau, x^{*}(\tau), u^{*}(\tau)) \right]\leq 0,
\]
meaning that 
\[
H\left(\tau, x^{*}(\tau), u^{*}(\tau), \lambda(\tau)\right) 
\geq H\left(\tau, x^{*}(\tau), v, \lambda(\tau)\right), 
\]
where $H= L(t, x, u) + \lambda \cdot f(t, x, u)$. 
Because $\tau$ is an arbitrary Lebesgue point 
of the control $u^{*}$ and $v$ is an arbitrary element of the set $\Omega$, 
it follows that the relation 
\[
H(t, x^{*}(t), u^{*}(t), \lambda(t))
= \underset{\omega \in \Omega} 
\max \, H\left(t, x^{*}(t), \omega, \lambda(t)\right)
\]
holds at all Lebesgue points, which ends the proof.
\end{proof}


\section{An Illustrative Example}
\label{sec:Ex}

As an example of application of our main result, 
let us consider the following distributed-order 
fractional optimal control problem: 
\begin{equation}
\label{exple}
\begin{gathered}
J[x(\cdot), u(\cdot)]= \int^{5}_{1} (1-3u(t))x(t)dt 
\longrightarrow \max,\\
^{C}\mathbb{D}^{\psi(\cdot)}_{{1}^{+}}x(t)
= u(t)x(t), \quad a.e. \quad t\in [1, 5],\\
u(t) \in [0, 2],\\
x(1)= x_a > 0,
\end{gathered}
\end{equation}
where the distribution function of order of differentiation
is given by
\[
\psi(\alpha)= \frac{\alpha}{3}.
\]

Let $u^{*}$ be an optimal control to problem \eqref{exple}.
Theorem~\ref{theo} give us a necessary optimality condition
that $u^{*}$ must satisfy. The Hamiltonian function associated 
with this problem is given by
\[
H(t, x, u, \lambda)= (1-3u)x+u x \lambda.
\]

From the maximality condition \eqref{maxcondition}, 
we know that $u^{*}(t)$ maximizes a.e. in $[0, 2]$ the mapping
\[
w\mapsto (1-3w)x^{*}(t) + wx^{*}(t)\lambda(t).
\]

Due to the positiveness of the initial condition ($x_a>0$) and 
the linearity of the distributed order derivative, we have that 
$x^{*}(t)>0$ for all $t\in [1, 5]$. Thus, the mapping to be
maximized can be reduced to 
\[
w\mapsto (\lambda(t)-3)w,
\]
and $u^{*}$ has the form
\begin{equation*}
u^{*}(t)=\left\{ 
\begin{array}{c}
2 \quad \text{ if } \, \lambda(t)>3,\\
0\quad \text{ if } \, \lambda(t)<3.
\end{array}
\right.
\end{equation*}

Now, it remains to determine the switching structure of the control through 
investigation of the adjoint boundary value problem given by 
\eqref{adj} and \eqref{trans}, that is, 
\begin{equation*}
\left\{ 
\begin{array}{l}
\mathbb{D}^{\psi(\cdot)}_{5^{-}} \lambda(t)
= 1 + (\lambda(t)-3)u^{*}(t),\\[3mm]
\mathbb{I}^{1-\psi(\cdot)}_{b^{-}}\lambda(5)=0.
\end{array}
\right.
\end{equation*}

Note that, because problem \eqref{exple} does not have a terminal phase constraint, 
the fractional transversality condition \eqref{trans} is simplified to $\lambda(5)=0$. 
Moreover, since $\lambda(\cdot)$ is a continuous function, there is $\xi>0$ such that 
$u^{*}(t)=0$ for all $t\in [5-\xi, 5]$. With this, we have that 
$\mathbb{D}^{\psi(\cdot)}_{5^{-}} \lambda(t)= 1$, 
and it follows, by backward integration, that
\[
\lambda(t)=\frac{(5-t)^{\bar{\alpha}}}{m\Gamma(\bar{\alpha} +1)},
\]
where $\displaystyle{m= \int^1_0\frac{\alpha}{3}d\alpha= \frac{1}{2}}$
and $\bar{\alpha}\in [0, 1]$. Noting that, for 
$$
c:= 5-(3m\Gamma(\bar{\alpha}+1))^{\frac{1}{\bar{\alpha}}} \in [7/2, 5[,
$$ 
we get $\lambda(c)=3$, we conclude that 
\begin{equation*}
u^{*}(t)=\left\{ 
\begin{array}{c}
2 \quad \text{ if } \, 0\leq t < c,\\
0\quad \text{ if } \, c\leq t \leq 5.
\end{array}
\right.
\end{equation*}


\section{Conclusions}
\label{sec:conc}

Recent applications and experimental data-analysis studies
have shown the importance of systems with ``diffusing diffusivity''
in anomalous diffusion, modeled with fractional, standard Brownian 
motions and distributed-order operators
\cite{rv2:Golan,rv2:C8SM02096E,rv2:Korabel}.
The theory of the calculus of variations for distributed-order
fractional systems was initiated in 2018 by Almeida and Morgado
\cite{ricardo}, and it has been extended by the authors 
in 2020 to the more general framework of optimal control \cite{faical}. 
There, we established a weak Pontryagin Maximum Principle (PMP), 
under certain smoothness assumptions on the space of admissible functions,
where the controls are not subject to any pointwise constraint \cite{faical}. 
The objective of the present article was to state and prove 
a strong version of the PMP for distributed-order
fractional systems, valid for general non-linear dynamics and $L^\infty$
controls and, in contrast with References \cite{ricardo,faical},
without assuming that the controls take values on all the Euclidean space. 
Our statement is as general as possible, and it encompasses 
the distributed-order calculus of variations of \cite{ricardo}
and the weak PMP of \cite{faical} as particular cases. Moreover,
in the analysis of a strong version of PMP, we emphasized 
the use of needle-like variations to control perturbations, 
dealing with controls taking values on a closed set
in a much larger class of admissible 
functions than in References \cite{ricardo,faical}. 
Our approach began by proving results on continuity of 
solutions due to needlle-like variations, and then followed 
by a differentiability result on the state solutions 
with respect to perturbed trajectories. The statement and the proof 
of the Pontryagin Maximum Principle are rigorously given. 
Finally, the new necessary optimality conditions 
were illustrated by a simple example for which an analytical
solution could be found. To deal with real optimal control problems 
of Nature, which are impossible to solve analytically, it is important
to develop numerical methods based on the fractional distributed-order
Pontryagin maximum principle here obtained. This will be subject
of future research.

\vspace{6pt} 

\authorcontributions{The authors equally contributed 
to this paper, read and approved the final manuscript.
Formal analysis, F.N. and D.F.M.T.; 
Investigation, F.N. and D.F.M.T.; 
Writing---original draft, F.N. and D.F.M.T.; 
Writing---review \& editing, F.N. and D.F.M.T. 
All authors have read and agreed to the published 
version of the manuscript.}

\funding{This research was funded by the Portuguese 
Foundation for Science and Technology (FCT),
grant number UIDB/04106/2020 (CIDMA).
Nda\"{\i}rou was also supported by FCT 
through the PhD fellowship PD/BD/150273/2019.} 

\dataavailability{Not applicable.}

\acknowledgments{The authors are grateful to three anonymous
reviewers for several valuable comments, which helped 
them to improve the manuscript.}

\conflictsofinterest{The authors declare no conflict of interest.
The funders had no role in the design of the study; in the collection, 
analyses, or interpretation of data; in the writing of the manuscript, 
or in the decision to publish the~results.} 

\end{paracol}


\reftitle{References}


\end{document}